\documentclass{amsart}%
\usepackage{amsmath}
\usepackage{amsfonts}
\usepackage{amssymb}
\usepackage{graphicx}
\usepackage{color}%
\providecommand{\U}[1]{\protect\rule{.1in}{.1in}}

\theoremstyle{plain}

\numberwithin{equation}{section}
\begin{document}
\title[Mathematics of the Genome]{Mathematics of the Genome}
\author{Indika Rajapakse}
\address{\noindent University of Michigan}
\email{indikar@umich.edu}
\author{Steve Smale}
\address{City University of Hong Kong}
\email{smale@cityu.edu.hk}
\thanks{Communicated by Michael Shub }
\date{}

\begin{abstract}
This work gives a mathematical foundation for bifurcation from a stable
equilibrium in the genome. We construct idealized dynamics associated with the
genome. For this dynamics we investigate the two main bifurcations from a
stable equilibrium. Finally, we give a mathematical proofs of existence and
points of bifurcation for the Repressilator and the Toggle gene circuits.

\end{abstract}
\maketitle

\noindent\textbf{Introduction:\medskip}

\noindent Development of a single cell into a multicellular organism involves
a remarkable integration of gene expression, molecular signaling, and
environmental cues. Emergence of specific cell types depends on selective use
of genes, a phenomenon known as gene regulation, at critical times during
development [1]. \medskip

\noindent Our goal is to build a mathematical foundation for gene regulation
and its dynamics. We believe this could lead to new insights in biology. \ In
fact new understandings of biology is our main goal.\medskip

\noindent Gene regulation dynamics traditionally takes the form:
\begin{equation}
\frac{dx}{dt}=F\left(  x\right)  \tag{1}%
\end{equation}

\noindent Firstly we will formally define and derive $F$ from the biology of
regulatory units with an intermediate step of a labeled oriented graph
network, starting with nodes as genes. Besides determining $F$, we will give a
structure that can be used to exhibit the equilibria and the periodic
solutions of (1) and some properties relating these components of the
dynamics. While the ideas are known in biology, the mathematical development
is new. \medskip

\noindent Secondly elementary bifurcation theory, especially the pitchfork,
will be developed from the point of view of the first bifurcation from a
stable equilibrium. The pitchfork bifurcation will be derived from scratch,
since we have not found the present literature on this subject satisfactory
for genome dynamics.\medskip

\noindent A main motivation for our work consists of two important advances in
biology, Gardner et al. on the "Toggle" [22] and Elowitz and Leibler on the
"Repressilator" [24]. See Sections 3 and 4 respectively.\medskip

\noindent Two and three-gene networks will be given where the pitchfork and
Hopf bifurcations are exhibited, starting from Hill functions, giving the
precise points of bifurcation even when distinct degradation rates of the
proteins are allowed. Here we take a stability point of view, in contrast to
much work in mathematical biology that emphasizes a Boolean point of view or
use of fixed point theory. In a future account we expect to give general
models of emergence of a single cell type as well as cell differentiation and
cell memory. \medskip

\noindent In 1967 Kauffman and McCulloch [2] introduced boolean networks as
models for genetic regulation, work that was extended by Kauffman in 1969 [3].
Similar models for regulatory networks have been discussed by Glass [4].
Related criteria for dynamics in feedback loops have been worked out for gene
networks [5]. Hastings et al. prove the existence of a periodic solution for
one class of $n-$dimensional systems using the Brouwer \ Fixed Point Theorem,
hence without stability information [6]. Although rigorous mathematical
treatments for cyclic feedback systems have been introduced by Mallet-Paret
and Smith [7] and Gedeon and Mischaikow [8], these works do not address the
stability questions as based on condition number (B\"{u}rgisser and Cucker,
[9]).\bigskip


\medskip\bigskip\bigskip

\noindent\noindent\textbf{Section 1:} \textbf{The network and associated
dynamics}\textit{.} \medskip

\noindent The configuration space is the set of $n$ nonnegative real numbers,
representing gene expression levels (for example measured by \textit{mRNA} or
protein abundance) of the $n$ genes. Thus, this configuration space\ is
$\left(
\mathbb{R}
^{+}\right)  ^{n}$ with its elements $x=(x_{1},\ldots,x_{n}).$ In fact we may
write $X=\prod_{i=1}^{n}\left[  0,k_{i}\right]  $, where $x_{i}$ represents
the concentration of a protein \textbf{\ }(i.e. protein abundance) and $k_{i}$
is its maximum value. Traditionally the dynamics are described by an ordinary
differential equation of the form
\begin{equation}
\frac{dx_{i}\left(  t\right)  }{dt}=F_{i}(x_{1}(t),\ldots,x_{n}(t)),i=1,...,n
\tag{2}%
\end{equation}
or in brief%
\[
\frac{dx}{dt}=F(x),\mathbf{\ }F\mathbf{:}X\mathbf{\rightarrow}%
\mathbb{R}
^{n}.
\]
\medskip The program of this section is: we will construct a network for the
genome. Toward that end we will define an oriented graph and then label the
nodes of the graph to obtain the network. Then $F$ will be derived from this
network.\medskip\medskip

\noindent\textbf{The Graph:\medskip}

\noindent The nodes are of two (distinct) types. First) the genes $g\in G.$
Second) the inputs $I$ \ to be interpreted as signals, as in Alon [10]. There
is also a special subset $T$ of $G,$ called transcription factors.\medskip

\noindent Now we describe the oriented edges of the graph, $E\left(
i,j\right)  $\ from node $i$ to node $j.$ The first type of edge satisfies
$i\in I$ and $j\in T.$ The second type satisfies $i\in T$ and $j\in G$ and is
interpreted as $i$ has some control over $j,$ as\textbf{\ }when the protein
generated by gene $i$\ \ is a transcription factor and can bind the promoter
of gene\textbf{\ }$j$\textbf{.} From this oriented graph we have, as usual in
graph theory, the concepts of paths and cycles. A predecessor of node $j$ is a
node $i$ such that $E\left(  i,j\right)  $ is a (well-defined) edge. We will
assume that each node $i\in I$ has no predecessor and the other nodes have at
least one predecessor. Also assume that the graph is connected [11]. \medskip

\noindent The "core" is \ defined to be the union of the following: nodes
which belongs to some cycle of the graph, and nodes with a path leading to a
cycle. \textbf{\medskip}

\noindent

\noindent Next we describe a layered structure on the nodes of this graph. The
zero$^{th}$ layer is to be taken as the core together with the set of inputs
\textbf{(}Figure 1C\textbf{)}. Inductively say that the node is in the
$k^{th}$\textbf{\ }layer\textbf{\ }provided all its predecessors lie in the
previous layer and it doesn't. Figure 1A, B, and C illustrate the cases of
graphs with empty core, a core with one layer, and a core with two layers,
respectively.\medskip\medskip\ \medskip%
\begin{figure}[ptb]%
\centering
\includegraphics[
natheight=6.301900in,
natwidth=8.000400in,
height=3.1488in,
width=3.9894in
]%
{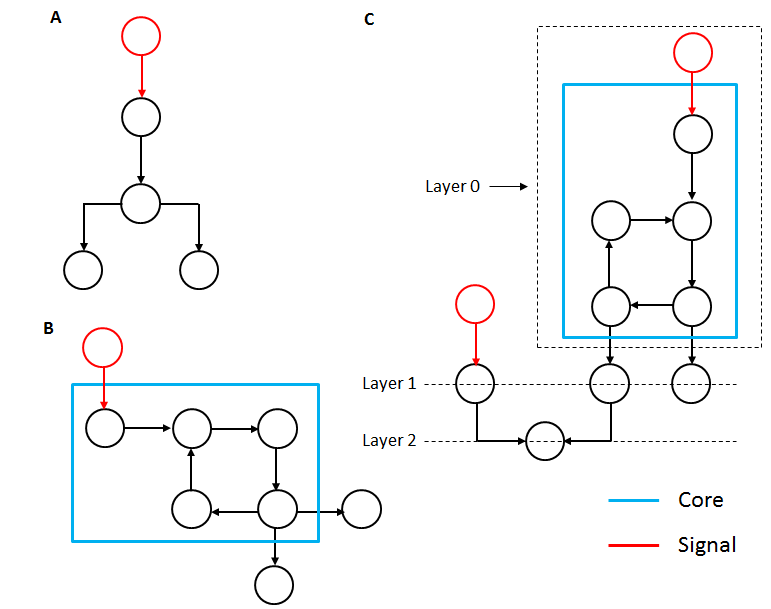}%
\caption{Graphs associated with the genome (A) Empty core, (B) Core with one
layer, (C) Core with two layers}%
\end{figure}

\noindent

\noindent\textbf{The Network:\medskip\ }

\noindent The next stage is to pass from the above graph to a "network" by
attaching some label to each node in $G.$ First consider a node whose
expression level is $y$ together with a single predecessor in $G$ whose
expression level is denoted by\textbf{\ }$x\in\left[  0,k\right]  .$ Inspired
by the theory of Hill functions [10, 12]\textbf{, }we let $H\mathbf{:}\left[
0,k\right]  \rightarrow%
\mathbb{R}
$ be given as well as a positive real number $b$. $H$\textbf{\ }is a map from
the concentration $x$\ to the rate of change of $y.$ Assume $H$\ is smooth,
nonnegative, and bounded by $bk.$ Here is the Hill type equation
\begin{equation}
\frac{dy}{dt}=H(x(t))-by(t),\text{ \ \ \ } \tag{3}%
\end{equation}
with the last term representing the degradation of $y,$ the level of
expression of the target node. Then we label the node that
expresses\textbf{\ }$y$\textbf{\ }with $\left(  H,b\right)  $\textbf{.
\ }Since our construction is so abstract it applies as well to the case that
the predecessor is an input node,\ or $x(t)$ is a signal. \ \medskip

\noindent

\noindent Next we extend the preceding Hill type function to the case of a
node $g$ with more than one predecessor. Define $P_{g}$ to be the set of
predecessors of $g\in G$ and let $W_{g}=\Pi_{i\in P_{g}}X_{i}$ where
$X_{i}=\left[  0,k_{i}\right]  $ is the space of expression of gene $i.$ The
Hill type function is a now a map
\begin{equation}
H_{g}:W_{g}\mathbf{\rightarrow}%
\mathbb{R}
,\text{ }\frac{dy}{dt}=H_{g}(x(t))-by(t). \tag{4}%
\end{equation}
For each gene $g,$ we label the node $g$\ with $\left(  H_{g},b\right)  .$
$H_{g}$ may be interpreted as aggregating the predecessor levels of $g$\ in
some way\textbf{.} We have completed the description of the network attached
to the genome.\medskip\ \ \ \medskip\noindent

\noindent\textbf{Network and the Dynamics:\medskip}

\noindent We can see this important fact: the network described by our graph
and the nodes just labeled (together with the inputs) determines the dynamics
of the genome. One just uses\textbf{\ }(4) to define the $F_{i}$ of $(1)$
where $F_{i}$ depends then only on the predecessors of gene $i.$
\ \textbf{\noindent}We identify gene $g$\ with gene $i$\ above and
below.\medskip\medskip

\noindent\textbf{Boundary Conditions:\medskip}

\noindent At a boundary point of $X,$ if $y_{i}=k_{i},$ some\textit{ }$i$,
then $\frac{dy_{i}}{dt}>0$ is a condition which is not compatible with our
interpretation of the biology\textbf{.} One cannot increase the concentration
if one is already at the maximum. This suggests that at such a point, $H_{g}$
$<bk_{g},$ which implies by Equation 4 $\frac{dy_{i}}{dt}<0.$ Now assume the
following hypothesis.\medskip

\noindent For each gene $g,$
\begin{equation}
\text{Function }H_{g}(x)\text{ in Equation 3 is bounded by }bk_{g}\text{ as a
function of }x. \tag{5a}%
\end{equation}

\noindent Moreover supplement 5a with%

\begin{equation}
H_{g}\left(  x\right)  >0,\text{ for }x\in W_{g}. \tag{5b}%
\end{equation}

\noindent Inequality 5b is important for a boundary point of $x$\ satisfying
$x_{i}=0$\ for some\textbf{\ }$i$\textbf{.} Then $\frac{dy_{i}}{dt}>0.$ \ The
boundary conditions imply the region $X$ is invariant under the forward
dynamics and in particular there are no equilibria on the boundary of\textbf{
}$X$.\bigskip

\noindent\textbf{Specific Hill Functions:\medskip}

\noindent Differential equations for gene control typically involve the Hill
functions, defined by either
\begin{equation}
H\left(  x\right)  =\frac{\beta x^{p}}{\left(  K_{1}\right)  ^{p}+x^{p}%
},\ \text{or }G(x)=\frac{\gamma}{1+\left(  \frac{x}{K_{2}}\right)  ^{q}},
\tag{6a}%
\end{equation}
where $x$\ now represents the concentration of a transcription factor in a
cell, and $H$\ represents activation, while $G$\ represents repression \ as
Alon [10]. The $p$ and $q$ might be considered as a function of a bifurcation
parameter $\mu$. \medskip\ 

\noindent In general $i$ is gene with a single predecessor $j$, $F_{i}\left(
x_{i}\right)  =H\left(  x_{j}\right)  -a_{i}x_{i}$ or $F_{i}\left(
x_{i}\right)  =G(x_{j})-a_{i}x_{i}$ as in Equation (2) \ and Equation (6a).
For all $i$, $H\left(  x_{i}\right)  =\frac{\beta_{i}x_{i}^{p}}{\left(
K_{1i}\right)  ^{p}+x_{i}^{p}}$ and similarly for $G(x_{i}).\medskip$

\noindent For both activation and repression at $p_{i}=q_{i}=0,$ $i=1,...,n,$
equilibrium equations simplify to $x_{i}=\frac{\beta_{i}}{2a_{i}}.$ \medskip
The equilibrium is globally stable since the eigenvalues of the Jacobian
matrix at the equilibrium are negative $\left(  -a_{i},i=1,...,n\right)  ;$
see Equation 7.\medskip

\noindent\noindent This equilibrium, $\left(  x\right)  _{i}=\left(
\frac{\beta_{i}}{2a_{i}}\right)  ,$ is the starting point of a path $x\left(
\mu\right)  $ of equilibria parameterized by a bifurcation parameter $\mu
<\mu_{0}.$ \medskip\ 

\noindent\medskip It is natural to ask that the equilibrium lies in the
feasible set according to the condition in 5a.\medskip

\noindent\noindent\noindent Now we can see the equilibrium discussed in
Equation 4a is feasible provided \
\begin{equation}
\frac{\beta_{i}}{2a_{i}}<k_{i}. \tag{6b}%
\end{equation}
This will be true if $k_{i}$ is big enough. \medskip

\noindent The aggregation of multipredecessor Hill functions can be
accommodated.\medskip

\noindent Sometimes we achieve greater generality by not using a specific form
of Hill functions. Moreover the mathematics are simpler.\medskip\medskip

\noindent\noindent\noindent\textbf{Induced equilibrium}: \medskip

\noindent Consider an equilibrium $x$ in the core dynamics of the network. It
has the form\textbf{\ }%
\[
x_{i}\left(  t\right)  =C_{i},\text{ a constant, for each }i\text{ belonging
to the core.}%
\]
Let gene $g$\ be in the first layer of the core. From Equation 4, $H_{g}$ at
the equilibrium in $W_{g}$, is therefore a constant say $C_{g}\,.$ Then, for
any expression $y(t)$ for gene $g$, \ $\frac{dy}{dt}=0$\ if and only if
$y(t)=\frac{C_{g}}{b}.$\ \ \ Now we define $y(t)\equiv$\ $\frac{C_{g}}{b}$. Do
the same for any gene $g$ in the first layer. By induction we do the same for
the whole genome. If $\ x$\ is a stable equilibrium the same will be true for
the corresponding $\ y$\ in the whole genome.\medskip\ This follows from the
Equation 4, since the derivative is $-b.$

\noindent\noindent We may say that an equilibrium solution in the core induces
a natural equilibrium in the set of all layers generated by the core in the
whole genome. The induced equilibrium is stable if the core equilibrium is
stable.\medskip\medskip\ 

\noindent\textbf{Remark 1}: We could replace the core by the zero$^{th}$ layer
in the above to obtain an equilibrium in the whole genome.\medskip
\smallskip\medskip

\noindent\textbf{Induced periodic solution:\medskip}

\noindent Suppose given a periodic solution $x(t)$ in the core with period
$T$. $\ $Consider a gene $g$ in the first layer. Project the periodic solution
$x(t)$ to $W_{g}.$ There exists a unique solution for $y\left(  t\right)  $ of
Equation 4 given an initial condition $y\left(  0\right)  $.\medskip

\noindent Rewrite Equation (4) as
\[
\frac{dy(t)}{dt}+by(t)=H_{g}\left(  x\left(  t\right)  \right)
\]
$\allowbreak$Multiply each side by $e^{bt}$ to yield%
\begin{align*}
\frac{d\left(  e^{bt}y\left(  t\right)  \right)  }{dt}  &  =e^{bt}H_{g}\left(
x\left(  t\right)  \right)  \text{.}\\
\text{for }s  &  >0,\text{ \ \ \ }e^{bs}y\left(  s\right)  =\int_{0}^{s}%
e^{bt}H_{g}\left(  x\left(  t\right)  \right)  dt+y\left(  0\right)
\end{align*}
for $s=T$ we have%
\[
e^{bT}y\left(  T\right)  =\int_{0}^{T}e^{bt}H_{g}\left(  x\left(  t\right)
\right)  dt+y\left(  0\right)  .
\]
If we set \noindent$y\left(  0\right)  =\frac{\int_{0}^{T}e^{bt}H_{g}\left(
x\left(  t\right)  \right)  dt}{e^{bt}-1},$ then from the equation above we
obtain \noindent$y\left(  T\right)  =y\left(  0\right)  .$ Now consider time
$kT+u,0\leq u<t,k$ is a non-negative integer.
\begin{align*}
e^{b\left(  kT+u\right)  }y\left(  kT+u\right)   &  =\int_{0}^{kT+u}%
e^{bt}H_{g}\left(  x\left(  t\right)  \right)  dt+y\left(  0\right)  .\\
&  =\int_{kT}^{kT+u}e^{bt}H_{g}\left(  x\left(  t\right)  \right)
dt+\overset{k-1}{\underset{m=0}{\sum}}\int_{mT}^{\left(  m+1\right)  T}%
e^{bt}H_{g}\left(  x\left(  t\right)  \right)  dt+y\left(  0\right) \\
&  =e^{bkT}\int_{0}^{u}e^{bt}H_{g}\left(  x\left(  t\right)  \right)
dt+\left(  \overset{k-1}{\underset{m=0}{\sum}}e^{bmT}\right)  \int_{0}%
^{T}e^{bt}H_{g}\left(  x\left(  t\right)  \right)  dt+y\left(  0\right) \\
&  =e^{bkT}\int_{0}^{u}e^{bt}H_{g}\left(  x\left(  t\right)  \right)
dt+\left(  \overset{k-1}{\underset{m=0}{\sum}}e^{bmT}\right)  \left(
e^{bT}-1\right)  y\left(  0\right)  +y\left(  0\right) \\
&  =e^{bkT}\left(  \int_{0}^{u}H_{g}\left(  x\left(  t\right)  \right)
dt+y\left(  0\right)  \right) \\
&  =e^{bkT}e^{bu}y\left(  u\right)  .
\end{align*}
Compare with the left side to get $y\left(  kT+u\right)  =y(u).$ Therefore,
oscillations in the core induce oscillations in the genome.

\bigskip

\noindent\noindent\textbf{Hardwiring:\medskip}

\noindent The network as described puts an oriented edge between genes $g$ and
$h$ if the protein product of gene $g$ \textit{can bind}, not will bind, to
the promoter of gene $h$ (in a broad sense) to promote the transcription of
$h$. Gene $g$ will bind only in some cell types, at certain stages of
development. It can happen that gene $g$ as a transcription factor may be
silenced. In that case gene $g$ can be removed from the network together with
its edges. For an example this phenomena can occur through chromatin
accessibility [13]. The universal gene network described in this section is
"hardwired" or is universal to gene interactions, just as the genes present in
the genome are invariant across all cell types.\textit{ \medskip
\smallskip\medskip}

\noindent\textbf{A stable equilibrium as a cell type: \ \medskip}

\noindent Gene regulatory networks can maintain characteristic patterns of
gene expression, thus defining the identity of a specific cell type. We
propose identification of a stable equilibrium (Definition 1)\textbf{ }for our
genome dynamics with a cell type. Both notions can be associated with unique
distributions (or proportions) of proteins. Thus, we consider a cell type that
is maintaining the gene expression pattern to be in a stable equilibrium. \ We
suggest that the emergence of a new cell type from this original cell type,
through differentiation or reprogramming, perhaps as a\ result of bifurcation,
is a departure from a stable equilibrium.\textbf{\ \medskip}

\noindent This phenomenon needs to be modified a bit for a stable oscillator
defining a cell type.

\bigskip

\noindent\textbf{Section 2:} \textbf{Bifurcations from a stable equilibrium in
}$n\mathbf{-}$\textbf{dimensional dynamical systems\medskip.}

\noindent We associate a pitchfork bifurcation with the phenomena of cell
differentiation or reprogramming in a number of situations based on analysis
of genomic data, but it is a challenging problem both in cell biology and
mathematics to solidify this process. The Hopf bifurcation is related to the
development of various oscillations in the genome that are often associated
with circadian rhythm. Vast genomic data exist for both cell differentiation
and circadian rhythms, and our motivation is to use these results to support a
more mathematical foundation for these processes. This in turn could guide a
better understanding of the genome and reprogramming. The rest of the section
is written towards developing such a mathematical\ foundation.\medskip

\noindent Our main finding is that there are exactly two generic bifurcations
as a first bifurcation from a stable equilibrium. These are (1) the
(supercritical) pitchfork bifurcation which produces two stable equilibria and
a separating saddle point and (2) the (supercritical) Hopf bifurcation which
produces a stable periodic solution.\medskip\noindent

\noindent Recall that the Jacobian matrix $J$\ of partial derivatives at an
equilibrium is: \textbf{\ }%
\begin{equation}
J=\left(
\begin{array}
[c]{ccc}%
\frac{\partial F_{1}}{\partial x_{1}} & ... & \frac{\partial F_{1}}{\partial
x_{n}}\\
\vdots & \frac{\partial F_{i}}{\partial x_{i}} & \vdots\\
\frac{\partial F_{n}}{\partial x_{1}} & ... & \frac{\partial F_{n}}{\partial
x_{n}}%
\end{array}
\right)  . \tag{7}%
\end{equation}
\noindent where $F=\left(  F_{1},...,F_{n}\right)  .\medskip$

\noindent\textbf{Definition 1: Stable equilibrium}. If\textbf{\ }$x_{0}$ is a
stable equilibrium (of Equation 2 for example) ; 1) all trajectories that
start near $x_{0}$ approaches it as $t\rightarrow\infty.$ That is, $x\left(
t\right)  \rightarrow$ $x_{0}$ as $t\rightarrow\infty.$ 2) Every eigenvalue of
the Jacobian matrix at $x_{0}$ has negative real part. 1) is a consequence of
2. The basin $B$ $\left(  x_{0}\right)  $\textit{ }is a set of all points
which tends to $x_{0}$ when $t\rightarrow\infty.\medskip$

\noindent\textbf{Definition 2: Closed orbits} (periodic solutions). They are
solutions for which $x\left(  t+T\right)  =x\left(  t\right)  $ for all $t,$
for some $T>0.$ A periodic solution is stable provided every nearby orbit
tends to it asymptotically as $t\rightarrow\infty.$ See [14] for eigenvalue
criteria for stability of closed orbits, as for Definition 1 as well.\medskip

\noindent These objects are fundamental in the subject of dynamics in
biology.\medskip\medskip\ 

\noindent\textbf{Definition 3}: \textbf{First bifurcation from a stable
equilibrium}. \ This is associated to a family $F_{\mu}$ with bifurcation
parameter $\mu\in(-e,e)$ describing a dynamic $\frac{dx}{dt}=F_{\mu}\left(
x\right)  .$ Here $x$ belongs to a domain $X$ of\textbf{\ }$%
\mathbb{R}
^{n}$ as in Section 1$,\ $and $F_{0}\left(  x\right)  =F\left(  x\right)  .$
We suppose that the dynamics of $F_{\mu}$ is that of a stable equilibrium
$x\left(  \mu\right)  ,$ basin $B_{\mu}$ for $\mu<\mu_{0},$ $\mu_{0}=0$ and
bifurcation at $\mu=0.\smallskip$ Let $J_{\mu}=$ $J\left(  x\left(
\mu\right)  \right)  $ as in Equation 7.\noindent\ We incorporate the
following hypothesis into definition 3. \ Suppose that the equilibrium does
not "leave it's basin" in the sense that there is a neighborhood $N$ of
$x_{0}$, such that $N$ is contained in $B_{\mu}$ for all $\mu<0.$
\bigskip\medskip

\noindent Note that the condition in the last statement of \ definition 3
implies by a uniform continuity argument, that even at $\mu=0$, $\ x_{0}$ is a
sink in the sense that $x\left(  t\right)  \rightarrow x_{0},$ if the initial
point belongs to $N$ for the dynamics of $F_{\mu},$ $\mu=0$ with basin
$B_{\mu_{0}}.\medskip$ It follows that in this space $N$ only equilibria at
$\mu$ is $x_{\mu}.\medskip$ This is our meaning of "first" in definition
3.\medskip

\noindent The eigenvalues of $J_{\mu}$\ at the stable equilibrium and before
the bifurcation value of $\mu$\ all have negative real part, either real, or
in complex conjugate pairs. At the bifurcation, one has either a single real
eigenvalue becoming zero and then positive with the pitchfork, or else a
complex conjugate pair of distinct eigenvalues with real parts zero to become
positive after the bifurcation and then the Hopf oscillation. A pitchfork
bifurcation converts a stable equilibrium into two stable equilibria. The Hopf
bifurcation converts a stable equilibrium into a stable periodic
solution.\medskip\smallskip

\noindent\noindent\textbf{Theorem 1}: Consider a first bifurcation from a
stable equilibrium with its corresponding notation (definition 3). Let
$D_{\mu}$ be the determinant of $J_{_{\mu}}$ at $x\left(  \mu\right)  .$ If
$D_{\mu_{0}}=0$ then generically the dynamics undergoes a pitchfork
bifurcation at $\mu_{0}.\medskip$

\noindent The "generically" for the case of one-variable is described with
$C_{3}\neq0$ in Lemma 1. \medskip\noindent Our use of the center manifold
theorem extends the notion to $n-$variables. \medskip

\noindent We give the proof of the Pitchfork bifurcation theorem, Theorem 1,
including the description of the bifurcation [15]. \noindent\medskip

\noindent\noindent\noindent We give the proof for one variable first. Then
$\frac{d\left(  F_{\mu}\left(  0\right)  \right)  }{dx}=\lambda$, the first
and only eigenvalue$.$ Let
\begin{equation}
F_{\mu_{0}}\left(  x\right)  =G\left(  x\right)  , \tag{8}%
\end{equation}
so that $G$ is independent of $\mu$ and $x_{0}$ taken to be zero$.\medskip$

\noindent\noindent\textbf{\noindent Lemma 1}: Let $C_{3}\neq0.$ The second
derivative of $F_{\mu_{0}}$ is zero$.\medskip$

\noindent\textbf{Proof:} \ Expand $F_{\mu_{0}}$ about zero in Equation 8 by
Taylor series.
\[
F_{\mu_{0}}\left(  x\right)  =C_{2}x^{2}+C_{3}x^{3}+R_{3}\left(  x\right)
\]
where $C_{2}=\frac{\text{D}^{2}}{2!}G\left(  0\right)  ,C_{3}=\frac
{\text{D}^{3}}{3!}G\left(  0\right)  ,$ where D $=$ $\frac{d}{dx}.$ $R_{3}$ is
the remainder and $R_{3}\left(  x\right)  \leq M\frac{\left\vert x\right\vert
^{4}}{4!},$ where D$^{4}\left(  x\right)  \leq M$ for all $\left\vert
x\right\vert \leq b,$ for appropriate $b.$ At the equilibrium $F_{\mu_{0}%
}\left(  x\right)  =0$, thus%
\[
0=C_{2}x+C_{3}x^{2}+R_{3}\left(  x\right)
\]
\noindent\noindent Zero is an equilibria of multiplicity two. A third
equilibrium is at $x=-\frac{C_{2}+R_{3}(x)}{C_{3}}.$ \medskip

\noindent There are two cases for which $C_{2}$ could be non-zero$.$ Take the
Taylor series remainder $\smallskip\medskip R_{3}\left(  x\right)  $ so that
$\left\vert R_{3}(x)\right\vert <\left\vert \frac{c_{2}}{2}\right\vert .$

\noindent\noindent1. $x=-\frac{C_{2}+R_{3}(x)}{C_{3}}$ belongs to the basin
$U=B_{0}.$ This contradicts the fact that $U$ is a basin of the
dynamics.\smallskip

\noindent\noindent\noindent2. $x=-\frac{C_{2}+R_{3}(x)}{C_{3}}$ does not
belong to the basin $U=B_{0}.$ This contradicts the Poincar\'{e}--Hopf index
theorem applied to the basin $U$ (see Section 5).\medskip

\noindent\noindent\noindent\noindent Therefore, the second derivative of $F$
at the equilibrium must be zero. Lemma 1 is proved.\medskip

\noindent We have supposed that $C_{3}\neq0$ by the genericity hypothesis and
then $C_{3}<0$ follows from the\ Taylor's formula together with the
bifurcation setting.\medskip

\noindent\noindent Let $J=J_{\mu}$ as in (7) \ with $\mu=0$ and $x=0.$
Proceeding to the general case, then by the genericity hypothesis the
eigenvalues $\lambda_{i}$ of $J$ could be assumed to be distinct. \ Then $J$
may be diagonalized with $\lambda_{1}=0$ and the real part of the $\lambda
_{i}$ is less than zero for $i>1.$ Perhaps this make the general case more
apparent.\medskip

\noindent The center manifold theory applies to give the reduction to
dimension one; the center manifold corresponds to the space of eigenvector of
eigenvalue zero\textbf{ }[16, 17]. \noindent\medskip

\noindent This finishes the proof of Theorem 1.\bigskip\ \noindent

\noindent Now \noindent consider the systems in 6a and 6b with all the
conditions and hypothesis as discussed. Also consider the bifurcation path of
equilibria parameterized by a bifurcation parameter $\mu<\mu_{0},$ starting at
the initial point described.$\medskip$

\noindent\textbf{Corollary of Theorem 1: }The first bifurcation as above is a
pitchfork, if at the first bifurcation, det$(J_{\mu})=0.\medskip$

\noindent\textbf{Proof of Corollary}: Note that the uniformity condition in
definition 3 (existence of $N$) follows from the fact that $x\left(
\mu\right)  $ is a first bifurcation. We may suppose that\textbf{ }$x\left(
\mu\right)  $\ is contained in $X$\textbf{. \ }This is because\textbf{
}$x\left(  \mu_{0}\right)  $ is not in the boundary of $X$\textbf{;
}even\textbf{ }$F\left(  \mu\right)  $ is never zero on the boundary. \medskip
This proves the Corollary, since the definition 3 is incoporated in Theorem
1.\medskip\ 

\noindent\noindent Previous treatments of the pitchfork assume a symmetry of
the function in Theorem 1a or else vanishing of its second derivative.
\medskip\noindent\noindent\ 

\noindent We suppose the following $F\left(  \mu\right)  :X\rightarrow%
\mathbb{R}
^{n},$ $X$ a domain in $%
\mathbb{R}
^{n}.\medskip\medskip$

\noindent\noindent\textbf{Hopf Bifurcation Setting: \medskip}

\noindent Let $x\left(  \mu\right)  $\ be a stable equilibrium for $\mu$\ for
all $\mu<\mu_{0}.$\ Suppose that the equilibrium does not "leave it's basin"
in the sense that there is a neighborhood $N$ of $x_{0}$, such that $N$ is
contained in $B_{\mu}$ for all\textbf{ }$\mu<0.$ Suppose that all eigenvalues
of Jacobian, $J$\ at $x\left(  \mu_{0}\right)  $\ have negative real parts
except one conjugate purely imaginary pair $\pm ib$ and these two eigenvalues
cross the imaginary axis when $\mu>\mu_{0}$. \noindent Then generically the
dynamics undergoes a bifurcation at $\left(  x_{0},\mu_{0}\right)
$\ resulting in a stable periodic solution. \ \bigskip This is called a Hopf bifurcation.

\noindent The "generically" of the periodic solutions involve two issues. 1)
$\frac{d\left(  \operatorname{Re}\lambda\left(  \mu\right)  \right)  }{d\mu
}\mid_{\left(  \mu=\mu_{0}\right)  }\neq0.$ 2) A certain third derivative is
not zero (this is first Lyapunov coefficient, of the dynamics $l_{1}%
\mid_{\left(  x_{0},\mu_{0}\right)  }$). See Mees and Chua [17], Marsden and
McCracken [18] , Guckenheimer and Holmes [15] for details.\medskip
\noindent\medskip

\noindent\textbf{Remark 2}: One can test the last condition using a numerical
continuation package, for example MATCONT [16].\noindent\ \medskip\noindent

\noindent Note that references 16, 17, 18, and 19 give us the justification
for Hopf Bifurcation Setting. \noindent While the original results by Hopf
[20] give conditions for the local existence of periodic orbits, a series of
papers by York and colleagues focused on the global continuation of periodic
orbits [21, 22]. Furthermore, in the paper by Mees and Chua, it was noted that
though the Hopf theorem only makes predictions locally, "experience tends to
confirm that predictions often remain qualitatively correct even when the
system is very far from bifurcation".\medskip

\noindent Normal forms for pitchfork and Hopf bifurcations \medskip\smallskip
are as follows.

\noindent\textbf{Pitchfork bifurcation:} \medskip

\noindent A normal form of the pitchfork bifurcation is%

\[
\frac{dx}{dt}=\mu x-x^{3},x\in%
\mathbb{R}
^{1},\ \ \mu\in%
\mathbb{R}
^{1}.
\]
$\mathbf{\ }$

\noindent This system has an equilibrium $x_{0}=0$\ for all $\mu.$\ This
equilibrium is stable for $\mu<0$\ and unstable for $\mu>0$\ ($\mu$\ is the
eigenvalue of this equilibrium). For $\mu>0,$\ there are two extra equilibria
branching from origin (namely, $x_{1,2}=\pm\sqrt{\mu})$\ which are stable.
\ This bifurcation is called a pitchfork bifurcation (see Figure 2). \medskip\ 

\noindent In $n-$dimensions we may write in addition, $\frac{dx_{i}}%
{dt}=-x_{i},i=2,...,n.\medskip\medskip$%

\begin{figure}[ptb]%
\centering
\includegraphics[
natheight=3.000000in,
natwidth=4.218600in,
height=2.7397in,
width=4.6458in
]%
{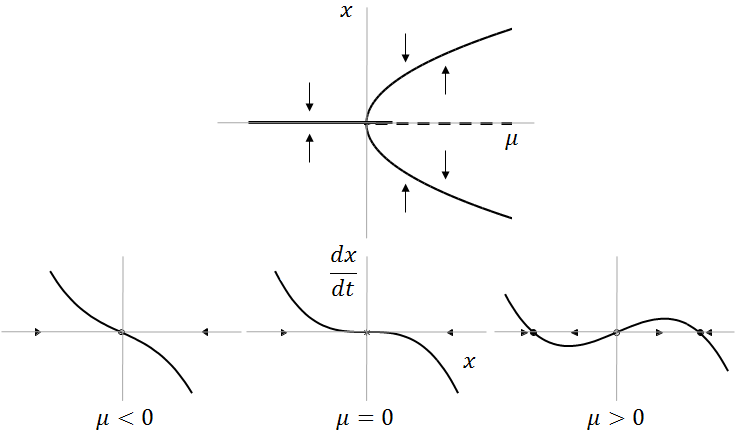}%
\caption{Pitchfork bifurcation (supercritical) diagram for $\frac{dx}{dt}=\mu
x-x^{3}.$ Top: The vertical axis is $x\left(  \mu\right)  ,$ the equilibria.
Solid lines represent stable, while the dotted line represents unstable.
Bottom: Three cases of $\mu.$ For $\mu<0$ there is a single stable fixed point
at $x=0.$ When $\mu=0$ there is still a single stable fixed point at $x=0$.
For $\mu>0,x=0$ is unstable, and flanked by two stable fixed points at
$x=\pm\sqrt{\mu}.$}%
\end{figure}

\noindent

\noindent\textbf{Hopf bifurcation: }\noindent\medskip

\noindent A main example of a Hopf\bigskip\noindent\ bifurcation is%
\begin{align*}
\frac{dx}{dt} &  =y-f_{\mu}\left(  x\right)  \\
\frac{dy}{dt} &  =-x.
\end{align*}

\noindent Here $f_{\mu}\left(  x\right)  =x^{3}-\mu x$ and $\mu$ lies in the
interval $\left[  -1,1\right]  .$ The only equilibrium in this system is
$x_{0}=0,y_{0}=0.$\ The Jacobian matrix $J$ at the equilibrium is%
\[
J=\left(
\begin{array}
[c]{cc}%
\mu & 1\\
-1 & 0
\end{array}
\right)  ,
\]

\noindent and the eigenvalues are $\lambda_{\pm}=\frac{1}{2}\left(  \mu
\pm\sqrt{\mu^{2}-4}\right)  .$ Thus the equilibrium is stable for $-1\leq
\mu<0$\ and has a unique stable periodic solution for $0<\mu<$ $c$, some $c$.
Figure 3 shows some phase portraits associated with this bifurcation.\medskip

\noindent In $n-$dimensions we may write in addition, $\frac{dx_{i}}%
{dt}=-x_{i},i=3,...,n,$ $x=x_{1},y=x_{2}.$%
\begin{figure}[ptb]%
\centering
\includegraphics[
natheight=2.052200in,
natwidth=8.895500in,
height=3.6089in,
width=3.8744in
]%
{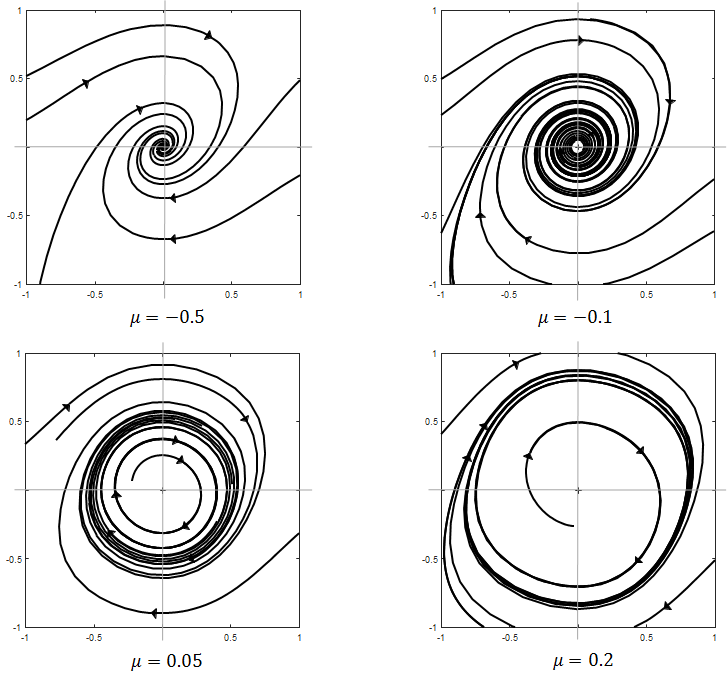}%
\caption{The Hopf bifurcation in the system \ $\frac{dx}{dt}=y-x^{3}+\mu
x,\frac{dy}{dt}=-x.$}%
\end{figure}

\bigskip

\noindent\textbf{Remark 3: }The literature discusses "supercritical" and
"subcritical" bifurcations. Because our bifurcations start with a stable
equilibrium, only supercritical bifurcations occur.\bigskip

\noindent In Section 3 and 4 we give examples of gene network dynamics, based
on the general theory of network dynamics developed in Section 1 and the
mathematics described in Section 2. \smallskip\medskip

\noindent\noindent\textbf{Section 3: Bifurcation analysis of two gene
networks\medskip.}\noindent

\noindent The Toggle, as designed and constructed by Gardner et al. [23], is a
network of two mutually inhibitory genes that acts as a switch by some
mechanism as a control switching from one basin to another.\medskip\ \ 

\noindent As stated in Gardner et al. "Here we present the construction of a
genetic toggle switch, a synthetic, bistable gene-regulatory network,
Escherichia coli and provide a simple theory that predicts the conditions
necessary for bistability." \ Ellner and Guckenheimer [24] developed some
mathematics of toggle of [23] in a way related to the development
here.\ \medskip\ 

\noindent We give a more a mathematical prospective of this phenomenon in the
setting of genome dynamics. \bigskip%
\begin{figure}[ptb]%
\centering
\includegraphics[
natheight=1.074100in,
natwidth=1.704500in,
height=0.6183in,
width=0.9651in
]%
{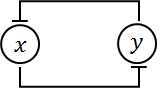}%
\caption{A schematic diagram of a two gene network with two repressors.}%
\end{figure}

\noindent\noindent Two gene network as in Figure (4), can be described by
\begin{align}
\frac{dx}{dt}  &  =F\left(  y\right)  -ax\tag{9}\\
\frac{dy}{dt}  &  =G\left(  x\right)  -by\text{\ }\nonumber
\end{align}
where $F\left(  y\right)  $ and $G\left(  x\right)  $ are Hill type functions
and these equations are special case of those in Equation 2.\medskip

\noindent The Jacobian matrix $J$ at an equilibrium $(x_{0},y_{0})$ is%
\[
J=\left(
\begin{array}
[c]{cc}%
-a & F^{^{\prime}}\left(  y_{0}\right) \\
G^{^{\prime}}\left(  x_{0}\right)  & -b
\end{array}
\right)  .
\]

\noindent In this case $F^{^{\prime}}<0,G^{^{\prime}}<0,$ since $F$ and $G$
are repressors$.$ Let $Q=F^{^{\prime}}G^{^{\prime}}$, thus $Q>0.$ \ The
characteristic equation of the matrix $J$ at $\left(  x_{0},y_{0}\right)  $ is
given by
\begin{equation}
\left(  \lambda+a\right)  \left(  \lambda+b\right)  =Q.\nonumber
\end{equation}
Eigenvalues at the equilibrium $\left(  x_{0},y_{0}\right)  $ are given by:
\begin{align}
&  \lambda_{1}=-\frac{1}{2}\left(  a+b\right)  +\frac{1}{2}\sqrt{\left(
a-b\right)  ^{2}+4Q}\tag{10}\\
&  \lambda_{2}=-\frac{1}{2}\left(  a+b\right)  -\frac{1}{2}\sqrt{\left(
a-b\right)  ^{2}+4Q}.\nonumber
\end{align}
\noindent

\noindent Note that $\det(J)=ab-Q$. \ So if $\det(J)=0,$ then $Q=ab$ and
$\lambda_{1}=0$ and $\lambda_{2}<0$.

\noindent\noindent\noindent Now we consider the particular setting by Gardner,
Cantor, and Collins of the toggle switch [23]. It consisted of two genes,
which they described as%

\begin{align}
\frac{dx}{dt}  &  =\frac{\alpha_{1}}{1+y^{m}}-x\tag{11}\\
\frac{dy}{dt}  &  =\frac{\alpha_{2}}{1+x^{n}}-y.\nonumber
\end{align}
These are particular case of Equation 9. If $m=n=0,$ the equilibrium is
$x=\frac{\alpha_{1}}{2},$ $y=\frac{\alpha_{2}}{2},$ and the eigenvalues of the
Jacobian are negative. Therefore, if $\alpha_{1}<2\max(x)$ and $<2\max(y),$
this equilibrium is feasible (as in Equation 6b). The system has a unique
global stable equilibrium. \ \medskip

\noindent Now we introduce the bifurcation setting. The bifurcation parameter
$\mu,$ for $\mu<\mu_{0},$ the system has a stable equilibrium. At $\mu_{0}$ we
assume a bifurcation. \ Take $m$ and $n$ to be bifurcating parameters, that is
$\left(  m,n\right)  =\left(  m_{\mu},n_{\mu}\right)  $, $\left(  m_{0}%
,n_{0}\right)  =(0,0),$ the initial value of the bifurcation parameter.
\medskip

\noindent If there is a first bifurcation with det$(J_{\mu})=0,$ by the
Corollary of Theorem 1, it is a pitchfork. Thus,\bigskip

\noindent\noindent\textbf{Theorem 2: }Suppose\textbf{ }$F_{\mu}$ is
represented by Equation 11 with\textbf{ }$\mu$ described as above. Let the
first bifurcation given by $\det(J_{\mu})=0.$ Then generically, there is a
pitchfork bifurcation at that $\mu$. \bigskip\ 

\noindent The existence of such a bifurcation is given in Appendix
1A.\medskip\medskip

\noindent\noindent\textbf{Section 4: Bifurcation analysis of three gene
networks\medskip.}\noindent

\noindent\noindent The gene networks discussed in this section include the
classical gene circuit, the Repressilator and an extension of the Toggle
discussed in the previous section [23, 25, 12]. The Repressilator is a model
for a biological circuit (an oscillating genetic circuit), as introduced by
Elowitz and Leibler [25], that consists of three genes that inhibit each other
in a cyclic way (Figure 5A). Extending this to a biological setting, they then
demonstrated that engineered Escherichia coli colonies exhibited the described
network oscillatory behavior. Three gene network related to a Toggle, defined
as a network of two mutually inhibitory genes and one excitatory gene (Figure
5B). \ Figure 5C is a network of two mutually excitatory genes and one
inhibitory gene, which produces behavior similar to a Repressilator.
\smallskip\bigskip\noindent\noindent\noindent%
\begin{figure}[ptb]%
\centering
\includegraphics[
natheight=2.603900in,
natwidth=6.364200in,
height=1.4598in,
width=3.5276in
]%
{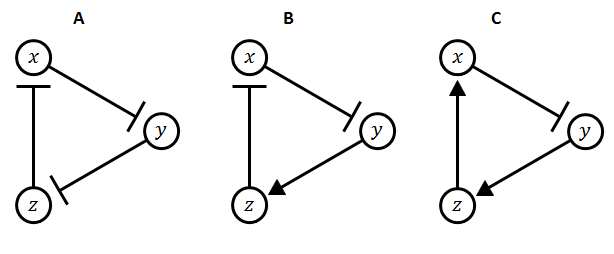}%
\caption{Gene circuits: A) A network with three inhibitory genes (the
Repressilator). B) A network with two inhibitory genes and one excitatory
gene. C) A network with two excitatory genes and one inhibitory gene. The
dynamics are the same as for network A). See below.}%
\end{figure}

\noindent Our three gene model is described by ordinary differential equations
of the form:
\begin{align}
\frac{dx}{dt}  &  =F\left(  z\right)  -ax\tag{12}\\
\frac{dy}{dt}  &  =G\left(  x\right)  -by\nonumber\\
\frac{dz}{dt}  &  =H\left(  y\right)  -cz\text{ \ }\nonumber
\end{align}

\noindent$x,y,z\in%
\mathbb{R}
^{+},$ and $a,b,c>0$ are degradation constants. \ Here $F,G,$\ and $H$\ are
Hill functions for the transcription factors shown in Figure 5.\bigskip\ \ 

\noindent We suppose that Equation 12 is of the type described in Equations 6a
and 6b. \ We are considering the "first bifurcation"\medskip.

\noindent The bifurcation setting is that described before Theorem 1. \ The
Jacobian matrix $J_{\mu}$ (as in Equation 7) at an equilibrium $(x_{0}%
,y_{0},z_{0})_{\mu}$ is%
\[
J_{\mu}=\left(
\begin{array}
[c]{ccc}%
-a & 0 & F^{^{\prime}}\left(  z_{0}\right) \\
G^{^{\prime}}\left(  x_{0}\right)  & -b & 0\\
0 & H^{^{\prime}}\left(  y_{0}\right)  & -c
\end{array}
\right)  .
\]
Let $Q=F^{^{\prime}}G^{^{\prime}}H^{^{\prime}}.$\ \ Det $\left(  J_{\mu
}\right)  =Q-abc.$ The characteristic equation is given by det$(J_{\mu
}-\lambda I)=0.$\ We expand the equation and obtain
\begin{equation}
\left(  \lambda+a\right)  \left(  \lambda+b\right)  \left(  \lambda+c\right)
=Q. \tag{13}%
\end{equation}

\noindent In this section, we similarly take $a=b=c=\alpha,$\ where
$\alpha=\frac{\left(  a+b+c\right)  }{3}$\ is the mean. \ In the subsequent
section, we eliminate this assumption.\medskip

\noindent\noindent The characteristic equation now is given by
\begin{equation}
\left(  \alpha+\lambda\right)  ^{3}=Q. \tag{14}%
\end{equation}

\noindent Assume that the system has a parameter $\mu$, that for $\mu<\mu_{0}%
$\ the equilibrium is stable and that it loses stability immediately after
$\mu=\mu_{0}.$\ \medskip Thus $Q=Q(\mu).\medskip$

\noindent\textbf{Remark 4}: The dichotomy of pitchfork and Hopf\ corresponds
to $Q>0,Q<0,$ respectively, where $Q$ is the value at $\mu_{0}$.\ \ \medskip

\noindent Let $m=$ largest real cube root of $\left\vert Q\right\vert
.$\medskip\noindent

\noindent\noindent\textbf{Case 1}: $Q<0.$ For example the Repressilator,
$F^{^{\prime}}<0,G^{^{\prime}}<0,H^{^{\prime}}<0.\medskip$ Another scenario
for Case 1 (see Figure 5C) is when\textbf{\ }$F^{^{\prime}}>0,G^{^{\prime}%
}>0,H^{^{\prime}}<0.$\noindent\medskip

\noindent Eigenvalues at the equilibrium $(x_{0},y_{0},z_{0})_{\mu}$ are given
by:
\[
\lambda_{1}=-\alpha-m,\ \lambda_{2}=-\alpha+\frac{1}{2}m+\frac{1}{2}i\sqrt
{3}m,\lambda_{3}=-\allowbreak\alpha+\frac{1}{2}m-\frac{1}{2}i\sqrt{3}m.
\]

\noindent Figure 6B and C illustrates the points $\lambda_{1}+\alpha,$
$\lambda_{2}+\alpha,$ and $\lambda_{3}+\alpha$\ corresponding to points in a
circle of radius one.\medskip\medskip

\noindent We again investigate the structure of the eigenvalues. When
$\frac{m}{2}<\alpha,$\ then $\lambda_{1}<0,$\ and the real part of
$\lambda_{2},\lambda_{3}<0.$\ When $\frac{m}{2}=\alpha,$ $\lambda_{1}<0,$\ and
$\lambda_{2},\lambda_{3}$\ are purely imaginary. When $\frac{m}{2}>\alpha
,$\ then $\lambda_{1}<0,$\ and the real part of $\lambda_{2},\lambda_{3}>0.$
These are sufficient for a generic Hopf bifurcation as in Section 2.
\ \bigskip

\noindent The existence of such a bifurcation is given in Appendix
1B.\medskip\medskip

\noindent\textbf{Case 2:} $Q>0$. For example, $F^{^{\prime}}<0,G^{^{\prime}%
}<0$ $,H^{^{\prime}}>0.$

\noindent Eigenvalues at the equilibrium $(x_{0},y_{0},z_{0})_{\mu}$
corresponding to the bifurcation value $\mu$ are given by:
\[
\lambda_{1}=-\alpha+m,\lambda_{2}=-\alpha-\frac{1}{2}m+\frac{1}{2}i\sqrt
{3}m,\lambda_{3}=-\alpha-\frac{1}{2}m-\frac{1}{2}i\sqrt{3}m.
\]

\noindent Now we investigate the structure of the eigenvalues carefully. When
$m<\alpha,$then $\lambda_{1}<0$\ and the real part of $\lambda_{2},\lambda
_{3}<0.$When $m=\alpha,$ $\lambda_{1}=0,$ real part of $\lambda_{2}%
,\lambda_{3}<0$\ $($since $-\alpha-\frac{1}{2}m<0).$\ If $m>\alpha$ then
$\lambda_{1}>0$\ and the real part of $\lambda_{2},\lambda_{3}<0.$\bigskip\ \ 

\noindent We are able to conclude the first bifurcation is a pitchfork by
Corollary of Theorem 1.\medskip\ 

\noindent The existence of such a bifurcation is given in Appendix
1C.\medskip\medskip\noindent

\noindent\noindent\textbf{Remark 5}: $\lambda_{1}+\alpha,\lambda_{2}%
+\alpha,\lambda_{3}+\alpha$\ are the three cube roots of $Q.$\ Points
$\lambda_{i}+\alpha$\ lie on a circle of radius $m$ .\ Figure 6A illustrates
these corresponding points in a circle of radius one.\medskip\ \noindent%
\begin{figure}[ptb]%
\centering
\includegraphics[
natheight=1.620700in,
natwidth=5.674000in,
height=0.9193in,
width=3.1479in
]%
{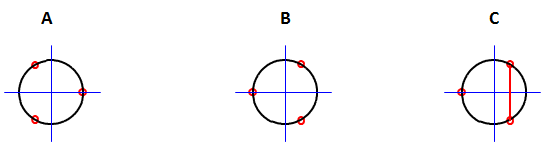}%
\caption{Two main bifurcations from a stable equilibrium for the three gene
network. A: $Q>0$ (pitchfork bifurcation), B: $Q<0$ (Hopf bifurcation). C:
illustration of the Hopf bifurcation point. The radius of the circle is one
and the marked points are the cube roots of 1 in A and of -1 in B.}%
\end{figure}

\noindent Here, we have given conditions for the two bifurcations for the 3
gene network.\medskip

\noindent\noindent We emphasize that this picture and these results depend on
the unrealistic hypothesis (as in the literature) that $a=b=c.$\ In the next
part we eliminate this hypothesis and exhibit the refined results in Theorem
3. \medskip\textbf{\ } \medskip\noindent\textbf{\ }

\noindent\textbf{Removal of the condition }$a=b=c.\smallskip\medskip$

\noindent\noindent Previous studies [23, 24, 12, 25, 26] have assumed that the
degradation constants are equal (one), but this hypothesis does not reflect
biological systems and besides can lead to ill-posed numerics. \medskip

\noindent Lets consider the setting of Equation 12 and the Jacobian $J_{\mu}.$
\smallskip

\noindent Then $Q=F^{^{\prime}}G^{^{\prime}}H^{^{\prime}}$ and note that the
determinant of $J=Q-abc$. The characteristic equation is given by (see
Equation 14)
\begin{equation}
(\lambda+a)(\lambda+b)(\lambda+c)=Q.\nonumber
\end{equation}
If we assume that there is a Hopf bifurcation at $\mu_{0}$, then a solution of
Equation 14 is a purely imaginary complex number and $Q<abc$. We may write
this solution in the form $\gamma\sqrt{-1},\gamma>0.$ Substituting this
solution into the characteristic equation, we get the following%
\begin{equation}
(\gamma\sqrt{-1}+a)(\gamma\sqrt{-1}+b)(\gamma\sqrt{-1}+c)=Q. \tag{15}%
\end{equation}
\bigskip\noindent\noindent\noindent Taking the real part of Equation 15 we get%

\begin{align}
-\left(  a+b+c\right)  \gamma^{2}+abc-Q  &  =0\nonumber\\
\gamma^{2}  &  =\frac{abc-Q}{\left(  a+b+c\right)  }>0 \tag{16}%
\end{align}

\noindent\noindent Taking the imaginary part of Equation 15 we get
\begin{align}
-i\gamma^{3}+iab\gamma+iac\gamma+ibc\gamma &  =0\nonumber\\
\gamma^{2}  &  =ab+ac+bc \tag{17}%
\end{align}
From Equations 16 and 17, we obtain%
\begin{align}
ab+ac+bc  &  =\frac{abc-Q}{\left(  a+b+c\right)  }\text{ }\tag{18}\\
\text{or }Q  &  =-2abc-\left(  a^{2}b+a^{2}c+ab^{2}+ac^{2}+b^{2}%
c+bc^{2}\right) \nonumber
\end{align}
Equation 18 is a necessary condition for Equation 12 to have a pair of pure
imaginary numbers as a solution.\medskip\smallskip

\noindent\textbf{Proposition 1}: Suppose there is a Hopf bifurcation at the
equilibrium $(x_{0},y_{0},z_{0})_{\mu=\mu_{0}}$. Then $abc-Q>0.$ Moreover,
\[
Q=-2abc-\left(  a^{2}b+a^{2}c+ab^{2}+ac^{2}+b^{2}c+bc^{2}\right)
\]
and the real eigenvalue is negative.\smallskip

\noindent\textbf{Proof}: \ The first two statements are in Equation 18. Let
$\lambda_{1}$ be the third eigenvalue; $\lambda_{1}=\frac{\det(J_{_{\mu_{0}}%
})}{\gamma^{2}},$ since the product of the eigenvalues is the determinant of
$J_{_{\mu_{0}}}$, and determinant is equal to $Q-abc.$ Thus, $\lambda
_{1}<0.\medskip$

\noindent On the other hand suppose we are given $a,b,c,$ and det$\left(
J_{_{\mu_{0}}}\right)  <0$ from the three gene problem (Equation 12). \ If
they satisfy Equation 18, then by taking $\gamma>0$ as in Equation 18 and
$\lambda_{1}$ as in Proposition 1, it follows that the eigenvalues are
$\pm\gamma\sqrt{-1}$ and $\lambda_{1}<0$. This satisfies the hypothesis for
the Hopf bifurcation. We suppose now that Equation 12 is of the type described
in Equations 6a and 6b. We have proved following theorem.\bigskip

\noindent\textbf{Theorem 3. }Consider Equation 15 and%
\begin{equation}
\left(  -2abc-\left(  a^{2}b+a^{2}c+ab^{2}+ac^{2}+b^{2}c+bc^{2}\right)
\right)  <Q<abc. \tag{19}%
\end{equation}

\noindent If $Q$ is in the above range in Equation 19, then there is no
bifurcation. The dynamic is that of a stable equilibrium. On the other hand,
if $Q$ equals the lower limit there is a Hopf bifurcation, and if $Q$ equals
the upper limit there is a Pitchfork bifurcation.\bigskip\medskip\ 

\noindent\textbf{Section 5: }$\mathbf{N-}$\textbf{dimensional theory\medskip.}

\noindent Here we assemble some brief remarks for the $n-$dimensional case.
This is some background \medskip for the problem of extending Theorem 3 to $n$
genes. \textbf{\medskip}

\noindent A network of order $n$. \ Consider the model described by ordinary
differential equations of the form:
\begin{align}
\frac{dx_{1}}{dt}  &  =F_{1}\left(  x_{n}\right)  -\alpha_{1}x_{1}\tag{20}\\
\frac{dx_{2}}{dt}  &  =F_{2}\left(  x_{1}\right)  -\alpha_{2}x_{2}\nonumber\\
&  \vdots\nonumber\\
\frac{dx_{n}}{dt}  &  =F_{n}\left(  x_{n-1}\right)  -\alpha_{n}x_{n}%
\text{\ }\nonumber
\end{align}

\noindent Let $x_{0}$ be an equilibrium and $J$ at $x_{0}.\smallskip$

\noindent Then let $Q=\prod F_{i}^{^{\prime}}$ and note that the determinant
of $J=-1^{n}\times Q-\prod\alpha_{i}.$ The characteristic equation for a
$n-$dimensional system can be written\textbf{\ }%
\begin{equation}
-1^{n}\overset{n}{\underset{i}{\prod}}\left(  \lambda+\alpha_{i}\right)
+Q=0,\ \tag{21}%
\end{equation}
where each $\alpha_{i}$ represents the absolute value of a diagonal term in
$J,$\ and is a degradation constant as in Equation 3. If we replace
$\alpha_{i}$ 's with the mean, that is $\alpha=$ $\frac{1}{n}\underset{i}{\sum
}\alpha_{i},$ then Equation 13 is simplified to the form.
\begin{equation}
\left(  \lambda+\alpha\right)  ^{n}=-1^{n+1}\times Q. \tag{22}%
\end{equation}
\noindent

\bigskip\noindent Using Equation 22, we can illustrate the two main
bifurcations, as shown in Figure 7 .%
\begin{figure}[ptb]%
\centering
\includegraphics[
natheight=4.337000in,
natwidth=5.674000in,
height=2.4128in,
width=3.1479in
]%
{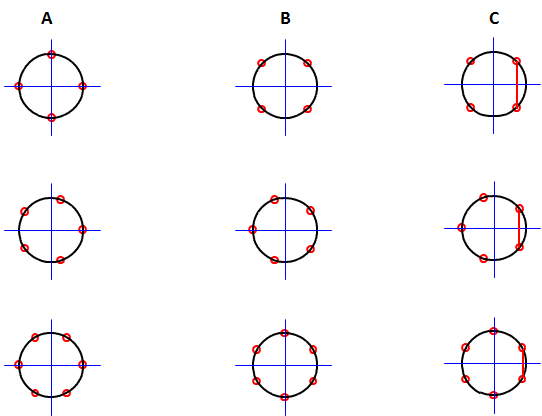}%
\caption{Two main bifurcations from a stable equilibrium. Column A: $Q>0$
(pitchfork bifurcation), Column B: $Q<0$ (Hopf bifurcation). Column C:
illustration of the Hopf bifurcation point. Rows show dimensions 4, 5 and 6
from top to bottom. The radius of the circle is one and the marked points are
the cube roots of 1 column A and of -1 in B and C.}%
\end{figure}
\qquad\medskip

\noindent\noindent\noindent\textbf{Poincar\'{e}-Hopf theorem}: \medskip

\noindent Recall our configuration space $X=\prod_{i=1}^{n}\left[
0,k_{i}\right]  $\ and the genome differential equation $\frac{dx}%
{dt}=F\left(  x\right)  $, $x\in X.$\ Then our boundary conditions guarantee
that no solution can leave $X.$\ Moreover, for each $x$ in the boundary of $X$
$,$ $F\left(  x\right)  $\ points strictly into $X.$\ More precisely, if we
ask that $F\left(  x\right)  $\ be strictly positive on $X$ then $F\left(
x\right)  $ points into the interior of $X$\ by Equation 6.\bigskip

\noindent\textbf{Theorem 4: }Suppose $\frac{dx}{dt}=F\left(  x\right)  $, $x$
belongs to $X$ and $F:X\rightarrow%
\mathbb{R}
^{n}.$ Suppose $X$ \ homeomorphic to open or closed ball [27].\bigskip%

\[
\sum_{F(x)=0,\text{ }x\in X}sign\left(  \det(J)\right)  \text{ }\left(
x\right)  =\left(  -1\right)  ^{n}%
\]
$J$\textbf{ }is the Jacobian of $F$ at $x.$\ \smallskip In particular
generically in the case $n$ is odd there is a odd number of equilibria.
\medskip

\noindent\noindent\noindent\textbf{Remark 6:} At a stable equilibrium in odd
dimensions, $\det(J)<0$ and in even dimensions $\det(J)>0.$\bigskip

\noindent\noindent A challenging problem is to demonstrate mathematically the
existence of a global stable oscillator in the Repressilator more generally
than what we have done here. Its existence has been shown computationally
[11]. This paper and the work of York and colleagues [21, 22] may provide a
framework for this problem.\medskip\smallskip

\noindent\textbf{Conclusion:\medskip}

\noindent In this work, we constructed a network for the genome, and defined
the concept of a core and layering in the network. Then we constructed a
dynamical system on gene expression levels from the network. For these
dynamics we investigated two main bifurcations from a stable equilibrium.
Last, we presented preliminary work toward understanding the bifurcations in
$n-$dimensions. We envision that this work will contribute to a deeper
understanding of the genome and its controllability.\bigskip

\noindent Note: Just prior to submitting the final revised version, we
received a related manuscript by Konstantin Mischaikow et al.\medskip

\noindent\textbf{Acknowledgments: }We would like to thank Yan Jun, Lu Zhang,
Lindsey Muir, Geoff Patterson, Xin Guo, Anthony Bloch, an anonmous referee and
especially Mike Shub for helpful discussions.\medskip\ 

\noindent We extend thanks to James Gimlett and Srikanta Kumar at Defense
Advanced Research Projects Agency for support and encouragement.\bigskip\

\noindent\textbf{Appendix 1A}: The existence of a pitchfork bifurcation in the
two gene case\medskip\noindent\ 

\noindent When $\alpha_{1}=\alpha_{2}=2,$ and $m=n>0,$ the system in Equation
11 can be written as%
\begin{align}
\frac{dx}{dt}  &  =\frac{2}{1+y^{m}}-x\tag{A1}\\
\frac{dy}{dt}  &  =\frac{2}{1+x^{m}}-y,\text{ }x,y\geq0.\nonumber
\end{align}

\noindent Proposition A: For our system given by A1, for $0\leq m\leq2,$ then
every equilibria must be $\left(  1,1\right)  .$ \medskip\noindent This can be
proved by \ numerics.\medskip

\noindent\noindent The Jacobian matrix $J$ at $(1,1)$%
\[
J=\left(
\begin{array}
[c]{cc}%
-1 & -m\frac{\alpha y^{m-1}}{\left(  y^{m}+1\right)  ^{2}}\\
-m\frac{\alpha x^{m-1}}{\left(  x^{m}+1\right)  ^{2}} & -1
\end{array}
\right)  =\left(
\begin{array}
[c]{cc}%
-1 & -\frac{1}{2}m\\
-\frac{1}{2}m & -1
\end{array}
\right)
\]
and the $det(J)=1-\frac{1}{4}m^{2}.$ Therefore, $0\leq m<2,$ $\ det(J)>0,$ and
$m=2$, $\ det(J)=0;$ this is a condition for a pitchfork bifurcation.\bigskip

\noindent\textbf{Appendix 1B: \ }The existence of a Hopf bifurcation in the
three gene case\medskip\noindent\ 

\noindent A simple case of the Repressilator can be written as\textbf{\ }%
\begin{align}
\frac{dx}{dt}  &  =\frac{\alpha}{1+z^{m}}-x\tag{B1}\\
\frac{dy}{dt}  &  =\frac{\alpha}{1+x^{m}}-y\nonumber\\
\frac{dz}{dt}  &  =\frac{\alpha}{1+y^{m}}-z,\text{ }x,y,z\geq0.\nonumber
\end{align}
The diagonal $x=y=z=s,$ is parameterized by $s.$ Then, the coordinates of an
equilibrium state are ($s,s,s),$ and the equilibrium path is described by%
\begin{equation}
s+s^{m+1}-\alpha=0. \tag{B2}%
\end{equation}
The Jacobian matrix $J$ at $(s,s,s)$ is%
\[
J=\left(
\begin{array}
[c]{ccc}%
-1 & 0 & -\frac{ms^{m}}{\left(  s^{m}+1\right)  }\\
-\frac{ms^{m}}{\left(  s^{m}+1\right)  } & -1 & 0\\
0 & -\frac{ms^{m}}{\left(  s^{m}+1\right)  } & -1
\end{array}
\right)  .
\]
The real part of the complex eigenvalues is, Re$(\lambda)=$ $-\frac{1}%
{2}\left(  \frac{-ms^{m}}{\left(  s^{m}+1\right)  }\right)  -1.$ The condition
for a Hopf bifurcation is thus
\[
-\frac{1}{2}\left(  \frac{-ms^{m}}{\left(  s^{m}+1\right)  }\right)  -1=0.
\]
This further simplifies to $s^{m}(2-m)=-2$ and for any $m>2,$ the following
equation determines a bifurcation value for $s.$
\[
s=\sqrt[m]{\frac{2}{m-2}}%
\]
Therefore this satisfies the conditions for a generic Hopf bifurcation. See
[25] \ for a related treatment.\medskip\medskip\medskip

\bigskip

\noindent\textbf{Appendix 1C}:\ The existence of a pitchfork bifurcation in
the three gene case\medskip\noindent\ 

\noindent As an example in Case 2 (Section 4), for the system described in
Equation 12 we take
\begin{align}
\frac{dx}{dt}  &  =\frac{2}{1+z^{m}}-x\tag{C1}\\
\frac{dy}{dt}  &  =\frac{2x^{m}}{1+x^{m}}-y\nonumber\\
\frac{dz}{dt}  &  =\frac{2}{1+y^{m}}-z.\nonumber
\end{align}
The equilibria of the system C1 are
\begin{equation}
x=\frac{2}{1+z^{m}},\text{ }y=\frac{2x^{m}}{1+x^{m}},\text{ and }z=\frac
{2}{1+y^{m}} \tag{C2}%
\end{equation}
\noindent Proposition C: For our system given by C1, $0\leq m\leq2,$ then all
equilibria $\left(  x,y,z\right)  $ must be $\left(  1,1,1\right)  .$ This can
be proved by numerics ( Stephen Lindsly). \medskip

\noindent The Jacobian matrix is%
\[
J=\left(
\begin{array}
[c]{ccc}%
-1 & 0 & -2m\frac{z^{m-1}}{\left(  z^{m}+1\right)  ^{2}}\\
2m\frac{x^{m-1}}{\left(  x^{m}+1\right)  ^{2}} & -1 & 0\\
0 & -2m\frac{y^{m-1}}{\left(  y^{m}+1\right)  ^{2}} & -1
\end{array}
\right)  =\left(
\begin{array}
[c]{ccc}%
-1 & 0 & -\frac{1}{2}m\\
\frac{1}{2}m & -1 & 0\\
0 & -\frac{1}{2}m & -1
\end{array}
\right)  .
\]
$\allowbreak$The $det(J)=$ $\frac{1}{8}m^{3}-1,$ which is zero exactly when
$m=2$. \noindent This satisfies the conditions for the existence of a
Pitchfork bifurcation.\bigskip\

\end{document}